\documentclass[11pt,leqno]{amsart}
\usepackage{amsmath,amssymb}

\usepackage{mathtools}
\usepackage{amsmath}
\usepackage{amssymb}
\usepackage{amsthm}
\usepackage{amscd}
\usepackage[msc-links,initials]{amsrefs}
\usepackage{bbm}
\usepackage{color}
\usepackage[english=american]{csquotes}
\usepackage{bm}
\usepackage{color}
\usepackage{enumerate}
\usepackage{mathrsfs}

\usepackage{hyperref}
\hypersetup{
	colorlinks=true,
	citecolor=blue,
	linkcolor=blue,
	filecolor=black,      
	urlcolor=blue,
}

\usepackage[shortlabels]{enumitem}

\newcommand{\mres}{\mathbin{\vrule height 1.6ex depth 0pt width
		0.1ex \vrule height 0.1ex depth 0pt width 0.5ex}}

\newcommand{\Acal}{\mathcal{A}}

\newcommand{\Ebf}{\mathbf{E}}
\newcommand{\Fbf}{\mathbf{F}}

\newcommand{\Mbf}{\mathbf{M}}
\newcommand{\Nbf}{\mathbf{N}}

\newcommand{\Rbf}{\mathbf{R}}

\newcommand{\Abb}{\mathbb{A}}

\newcommand{\Nbb}{\mathbb{N}}

\DeclareMathOperator{\im}{Im}

\DeclareMathOperator{\diverg}{div}

\DeclareMathOperator{\rank}{rank}

\DeclareMathOperator{\spn}{span}

\newcommand{\set}[2]{\left\{\, #1 \ \textup{{:}}\ #2 \,\right\}}

\newcommand{\dpr}[1]{\langle #1 \rangle}

\newcommand{\cl}[1]{\overline{#1}}

\newcommand{\R}{\mathbf{R}}

\newcommand{\loc}{\mathrm{loc}}
\newcommand{\sym}{\mathrm{sym}}

\newcommand{\embed}{\hookrightarrow}

\newcommand{\sbullet}{\begin{picture}(1,1)(-0.5,-2)\circle*{2}\end{picture}}
\newcommand{\frarg}{\,\sbullet\,}

\newcommand{\eps}{\epsilon}

\DeclareMathOperator{\Hom}{Hom}

\DeclareMathOperator{\Tan}{Tan}

\newtheoremstyle{thmlemcorr}{15pt}{15pt}{\itshape}{}{\itshape\bfseries}{.}{5pt}{{\thmname{#1}\thmnumber{ #2}\thmnote{ (#3)}}}
\newtheoremstyle{thmlemcorr*}{15pt}{15pt}{\itshape}{10pt}{\bfseries}{.}\newline{{\thmname{#1}\thmnumber{ #2}\thmnote{ (#3)}}}
\newtheoremstyle{defi}{15pt}{15pt}{}{}{\bfseries}{.}{5pt}{{\thmname{#1}\thmnumber{ #2}\thmnote{ (#3)}}}
\newtheoremstyle{remexample}{15pt}{15pt}{}{}{\itshape}{.}{5pt}{{\thmname{#1}\thmnumber{ #2}\thmnote{ (#3)}}}
\newtheoremstyle{ass}{10pt}{10pt}{}{}{\bfseries}{.}{}{{\thmname{#1}\thmnumber{ }\thmnote{ (#3)}}}

\theoremstyle{thmlemcorr}
\newtheorem{theorem}{Theorem}
\newtheorem{corollary}[theorem]{Corollary}
\newtheorem{proposition}[theorem]{Proposition}

\theoremstyle{thmlemcorr*}
\newtheorem{theorem*}{Theorem}
\newtheorem{lemma*}[theorem]{Lemma}
\newtheorem{corollary*}[theorem]{Corollary}
\newtheorem{proposition*}[theorem]{Proposition}
\newtheorem{problem*}[theorem]{Problem}
\newtheorem{conjecture*}[theorem]{Conjecture}

\theoremstyle{defi}

\theoremstyle{remexample}
\newtheorem{remark}{Remark}

\title[A Lebesgue--Lusin property]{A Lebesgue-Lusin property for linear operators of first and second order}
\author{Adolfo Arroyo-Rabasa}
\email{adolforabasa@gmail.com}
\address{Universit\'e catholique de Louvain, Belgium}

\begin{document}

%\title[A Lebesgue--Lusin property for first-order operators]{A Lebesgue--Lusin property \\ for gradients of first-order operators}
%
%\author{Adolfo \surname{Arroyo-Rabasa}}
%
%\address{First author address \email{adolforabasa@gmail.com}}

\begin{abstract}
We prove that for a homogeneous linear partial differential operator $\Acal$ of order $k \le 2$ and an integrable map $f$ taking values in the essential range of that operator, there exists a function $u$ of special bounded variation satisfying
	\[
		\Acal u(x)= f(x)  \qquad \text{almost everywhere}.
	\]
\noindent This extends a result of G. Alberti for gradients on $\R^N$. In particular, for $0 \le m < N$, it is shown that every integrable $m$-vector field is the absolutely continuous part of the boundary of a normal $(m+1)$-current.
\end{abstract}

%\keywords{bounded variation, Lusin property, first-order differential operator}
%
%\classification{49Q15; 35F05; 26B30}

\maketitle

	\section{Introduction}
Let $\Omega$ be an open subset of $\R^N$. We consider a general constant-coefficient system of linear partial differential equations acting on functions $u : \Omega \to E$, that is, a partial differential operator of the form 
\begin{equation}\label{eq:A}
\Acal u \, = \, \sum_{|\alpha| = k} A_\alpha \partial^\alpha u, \qquad A_\alpha \in   \Hom(E,F) \, ,
\end{equation} 
where $E,F$ are finite-dimensional $\R$-spaces. Here $\alpha = (\alpha_1,\dots,\alpha_N) \in (\Nbb_0)^N$ is a multi-index with modulus $|\alpha| = \alpha_1 + \cdots + \alpha_N$ and $\partial^\alpha = \partial_1^{\alpha_1} \circ \cdots \circ \partial_N^{\alpha_N}$ is the composition of distributional directional partial  derivatives.
The Fourier transform establishes a one-to-one correspondence between homogeneous operators and their associated principal symbol $\Abb: \Rbf^N \to \Hom({E,F})$, which in this context is given by the $k$-homogeneous tensor-valued polynomial
	\[
		\mathbb A^k(\xi) = \sum_{|\alpha| = k} A_\alpha \xi^\alpha, \qquad \xi^\alpha = \xi_1^{\alpha_1} \cdot \cdots \cdot \xi_N^{\alpha_N}, \quad \xi \in \R^N \, .
	\]

Suppose now that $u \in \mathscr D'(\Omega; E)$ is a distribution such that $\Acal u$ is a zero-order distribution, i.e.,  represented by an $F$-valued Radon measure. Let $g: \Omega \to F$ be the polar vector of $\Acal u$. If, additionally, $\Acal$ is a constant-rank operator (see, e.g.,~\cites{Murat,Raita,AS}),  the De Philippis--Rindler theorem~\cite[Theorem~1.1]{DR} establishes that $g$ is constrained to take values on the \emph{image cone} at \emph{singular points}. More precisely, if $U \subset \Omega$ is a $\mathscr L_N$-negligible Borel set, where $\mathscr L_N$ stands for the $N$-dimensional Lebesgue measure, then
\[
	g(x) \in \mathcal I_{\Acal} \coloneqq \bigcup_{|\zeta| = 1} \im \mathbb A^k(\zeta) \qquad \text{$|\Acal u|$ almost everywhere on $U$\,,}
\]
This property extends a classical result of Alberti~\cite{rank_one}, which says that if the distributional gradient $Du$ of a map $u: \R^N \to \R^M$ is represented by an $\R^{M \times N}$-valued Radon measure, then its polar vector must take values in the cone of rank-one matrices on $\mathscr L_N$-negligible sets: if $\mathscr L_N(U) = 0$, then
\[
	\rank\, \frac{Du}{|Du|}(x)  = 1 \qquad \text{$|Du|$ almost everywhere in $U$.}
\] 
In contrast with this restriction over the polar of $Du$ at singular points, Alberti showed~\cite[Theorem~3]{Alberti} that the absolutely continuous part of a gradient measure is fully unconstrained: if $f \in L^1(\R^N;\R^{M \times N})$, then there exists $u \in BV(\R^N;\R^M)$ such that \footnote{Alberti's result is written originally for functions $f: \Omega \to \R$. Applying it on each component yields the general case for $\R^M$-valued vector fields.}
\begin{equation}\label{eq:Alberti}
	Du = f \, \mathscr L_N + [u] \otimes \nu_u \, \mathscr H_{N-1} \mres J_u\,, \qquad {[u] \coloneqq u^+ - u^-\,},\end{equation}
	and
	\[
		 \|u\|_{BV} \le C \int_\Omega |f|\,.
	\]
Alberti (see~\cite{Alberti,Alberti2}) also established \emph{other} {Lusin}-type properties for gradients of arbitrary order. He showed that if $f : \Omega \to \R^{N^k}_\sym$ is continuous, then for any positive measure $\lambda$ on $\Omega$ and any smallness constant $\eps > 0$, one may find a compact set $K \subset \Omega$ and a function $u \in C^k(\Omega)$ such that 
\[
	f(x) = D^k u(x) \quad \forall x \in K\,, \qquad \lambda(\Omega \setminus K) \le \eps \lambda(\Omega)\,,
\] 
and satisfying $L^\infty$-estimates for $D^k u$ in terms of $f$. Building on these ideas, Francos~\cite{Greg} showed that for any given Borel $f : \Omega \to \R^{N^k}_\sym$ and $\sigma > 0$, there exists a function $g \in C^{k-1}(\Omega)$, $k$-times differentiable almost everywhere, satisfying
\[
	f(x) = D^k g(x) \quad \text{a.e. on $\Omega$} \; \quad \text{and} \; \quad \|D^k g\|_{L^\infty(\Omega)} \le \sigma. 
\]
Driven by applications to higher-order variational problems where derivatives give rise to surface energies, Fonseca et al.~\cite[Theorem 1.4]{Fonseca} established an analog of~\eqref{eq:Alberti} for the Hessian operator. More precisely, they showed that {if} $f : \Omega  \to \R^{N \times N}_\sym$, then there exists $u \in W^{1,1}(\Omega)$, with $\nabla u \in BV(\Omega;\R^{N})$, satisfying
	\[
		D^2 u = f \, \mathscr L_N \, + \, [\nabla u] \otimes \nu_u \, \mathscr H_{N-1} \mres J_{\nabla u} \; \quad \text{and} \; \quad \|u\|_{BH(\Omega)} \le C \|f\|_{L^1}\,.
	\]
Unlike the Lusin-type properties for higher-order gradients by Francos, this generalization to second-order derivatives is not followed by a straightforward iteration of Alberti's property~\eqref{eq:Alberti}. The reason for this drawback is the presence of the symmetry constraints of higher-order curl-free fields.

In light of this digression, the recent developments in the study of fine properties in $BV^\Acal$ spaces (see, e.g.,~\cites{AD,Slicing,Anna,Breit,Diening}), and the De Philippis--Rindler theorem, we are led to ask the following natural question: \emph{is the absolutely continuous part of an $\Acal$-gradient measure fully unconstrained?} Our main results (Theorems~\ref{thm:1} and~\ref{thm:2}) establish that, \emph{at least} when $\Acal$ is a an operator of order $k \le 2$, this is indeed the case; regardless of $\Acal$ satisfying the constant-rank condition. Our proof hinges on Proposition~\ref{prop}, where we show the previous question is equivalent to the validity of~\eqref{eq:Alberti} for Hessians of arbitrary order. This is the main reason why our result is restricted to the case $k = 1,2$. %Albeit we believe that the techniques in~\cite{Fonseca} are likely to be extended to the higher-order case, it is not our objective to pursue this goal here.

\section{Notation and results} 
In all that follows, $\Omega$ will denote an open subset of $\Rbf^N$. We will write $\mathscr L_N$ to denote the $N$ dimensional Lebesgue measure and $\mathscr H_{N-1}$ to denote the ($N-1$)-dimensional Hausdorff outer measure on $\Rbf^N$. The space $BV(\Omega;E)$, of $E$-valued vector fields with bounded variation on $\Omega$, consists of all integrable maps $u: \Omega \to E$ whose distributional gradient can be represented by a finite Radon measure taking values on $E \otimes \R^N$. For such functions, we shall write $\nabla u \in L^1(\Omega;E \otimes \R^N)$ to denote the absolutely continuous part of $Du$ with respect to $\mathscr L_N$.
It is well known (see, e.g.,~\cite{APF}) that such maps are Lebesgue continuous outside of a countably $\mathscr H_{N-1}$ rectifiable set $J_u \subset \Omega$, with orientation normal $\nu_u$, called the jump set of $u$. Moreover, the map $u$ has Lebesgue one-sided limits $u^+(x),u^-(x)$ for $\mathscr H_{N-1}$ almost every $x \in J_u$ with respect to the normal direction $\nu_u(x)$, and their difference
\[
[u](x) \coloneqq u^+(x) - u^-(x), \qquad x \in J_u,
\]
defines a $\mathscr H_{N-1}$-integrable map on $J_u$. 
By the Radon--{Nikodym} theorem, the gradient of a map $u \in BV(\Omega;E)$ can be decomposed as $Du = \nabla u \, \mathscr L_N + D^s u$ where $D^s u \perp \mathscr L_N$. The singular part $D^s u$ can be further decomposed into mutually singular measures as
\[
	D^s u = D^c u \, + D^j u\,,
\]
where (the \emph{Cantor part}) $D^c u$ is a measure vanishing on countable unions of $\mathscr H_{N-1}$-finite Borel sets, and the so-called \emph{jump part} 
\[
D^j u \coloneqq [u] \otimes \nu_u \, \mathscr H_{N-1} \mres J_u
\]
is the restriction of $Du$ to the set of jump discontinuities of $u$. The subspace $SBV(\Omega;E)$ of $E$-valued functions of special bounded variation consists of all functions $u \in BV(\Omega;E)$ with rectifiable singular gradient, i.e., 
\[
	Du = \nabla u\, \mathscr L_N + [u] \otimes \nu_u \, \mathscr H_{N-1} \mres J_u\,.
\]
	
 We also define
\begin{equation}\label{eq:F}
\Fbf_{\Acal} \coloneqq \cl{\set{\Acal u(x)}{{x\in \Omega}, \, u \in \mathscr D(\R^N;E)}^F} \, ,
\end{equation}
the \emph{essential range} of the operator $\Acal$. It is easy to see that every distribution $\Acal u$ takes values in this space, which is the smallest space with this property. 

With these considerations in mind, we can state our main results. We begin by stating the result for first-order operators:
	
	\begin{theorem}\label{thm:1} Let $\Acal$ be a first-order operator as in~\eqref{eq:A} and let  $f : \Omega \to \Fbf_{\Acal}$ be integrable. 
	Then there exists a map $u \in {SBV}(\Omega;E)$ satisfying 
	\begin{gather}
	\Acal u \; = \; f\,\mathscr L_N \; + \;  \Abb^1(\nu_u)[u] \, \mathscr H_{N-1} \mres J_u\,
	%|Du|(\Omega) \le C \|f\|_{L^1(\Omega)}
	\end{gather}
	and
	\[
		\int_\Omega (|u| + |\nabla u|) \, dx + 	\int_{J_u} |[u]| \, d\mathscr H_{N-1} \le C \int_\Omega |f|\, dx
	\]
\noindent for some constant $C$ that only depends on $\Acal$.
	\end{theorem}
	
		A similar statement holds for second-order operators:
	
\begin{theorem}\label{thm:2} Let $\Acal$ be a second-order operator as in~\eqref{eq:A} and let  $f : \Omega \to \Fbf_\Acal$ be integrable. 
	Then there exists a map $u \in W^{1,1}(\Omega;E)$, with $\nabla u \in SBV(\Omega;E \otimes \R^N)$, satisfying
	\begin{gather}
	\Acal u \; = \; f\,\mathscr L_N \; + \;  \Abb^2(\nu_{{\nabla u}})( \,[\nabla u] \nu_{\nabla u}\,) \, \mathscr H_{N-1} \mres J_{\nabla u}
	\end{gather}
 and 
\[
\int_\Omega (|u| + |\nabla u| + |\nabla^2 u|) \, dx + \int_{J_{\nabla u}} |[\nabla u]| \, \mathscr H_{N-1} \le C  \int_\Omega |f|\, dx
\]
\noindent where $C$ is a constant that only depends on $\Acal$.
	\end{theorem}

	\begin{remark}\label{rem:coefficients}The essential range $\Fbf_\Acal$ defined in~\eqref{eq:F}  coincides  with the space 	$\spn \set{\Abb^k(\xi)[e]}{\xi \in \R^N, e \in E} = \spn \, \{\mathcal I_{\Acal}\}$
of all $\Acal$-gradient amplitudes in Fourier space (for a proof see~\cite[Sec.~2.5]{ADR}). For example, if $\Acal = D^2$ is the Hessian, then 
$\Fbf_{D^2} = \spn \set{\xi \otimes \xi}{ \xi \in \R^N} \eqqcolon \R^{N \times N}_\sym$.
\end{remark}
Theorems~\ref{thm:1} and~\ref{thm:2} will follow directly from the equivalence result below (Proposition~\ref{prop}), Alberti's original result for gradients (in the first-order case) and its analog for Hessians by Fonseca, Leoni and Paroni (in the second-order case). Before stating Proposition~\ref{prop}, it will be convenient to introduce some basic notation for symmetric tensors and gradients of higher order. 
	
\subsection*{Symmetric tensors}Let $r \ge 0$ be an integer. For a vector $v \in \R^N$, we write $v^{\otimes ^r}$ to denote the tensor that results by taking the tensorial product of $v$ with itself $r$-times (with the convention $v^{\otimes^0} = 1$).  We consider the subspace
\begin{align*}
	\R^{N^r}_\sym & = \spn \set{v^{\otimes^r}}{v\in\R^N} %\subset  \underbrace{\R^N \otimes \cdots \otimes \R^N}_{\text{$r$-times}}
\end{align*}
consisting of all symmetric $r$th order tensors on $\R^N$.  In the following, we shall consider a contraction $\dpr{\frarg,\frarg}_r : (E \otimes \R^{N^r}_\sym \times \R^N) \to E$ associated with the canonical inner product $(\frarg , \frarg)$ on $\otimes^r \R^N$ as
\[
	\dpr{e \otimes V,v}_r = e (V, v^{\otimes^r}) \qquad \text{for all $e \in E$, $V\in \R^{N^r}_\sym$ and  $v\in\R^N$.}
\]
Notice that under our convention for $r=0$, we have $\dpr{e,\frarg}_0 = e$.

%For $P \in E \otimes \R^{N^k}_\sym$ and $v \in \R^N$, we define
%\[
%		P \bullet^r v =  \begin{cases}
%		\dpr{P, v^{\otimes^r} } & \text{if $r > 1$,}  \\
%		P & \text{if $r =0$,}
%		\end{cases}
%	\]
%where  
\subsection*{Jump densities for higher-order $BV^k$-spaces} Let $k \ge 2$ be an integer. We define the space $BV^k(\Omega)$ as the space of integrable functions $u: \Omega \to \R$ whose distributional $k$\textsuperscript{th} order gradient $D^ku$ can be represented by a Radon measure taking values on $\R^{N^k}_\sym$; the spaces $BV^k_\loc(\Omega)$, $BV^k(\Omega;E)$ are defined accordingly with this definition in the obvious manner. By classical elliptic regularity theory, there is a natural continuous embedding $BV^k_\loc(\Omega;E) \embed W^{k-1,1}_\loc(\Omega;E)$. Therefore, the distributional gradient $D^{k-1} u$ of a map $u$ in $BV^k_\loc(\Omega;E)$ can be represented by an integrable map $\nabla^{k-1}u : \Omega \to E \otimes \R^{N^{k-1}}_\sym$, that is,
\[
	D^{k-1}u = \nabla^{k-1} u \, \mathscr L_N\,.
\]
Moreover, in this case, the tensor field $w\coloneqq \nabla^{k-1} u$ has bounded variation, the jump set $J_w$ of $w$ is a countably $\mathscr H_{N-1}$ rectifiable set where $w$ has approximate one-sided limits $[w]\in E \otimes \R^{N^{k-1}}_\sym$, with respect to a fixed orientation $\nu_w$ of $J_w$, and the difference map $[w]$ is integrable on $J_w$. Since $D^k u$ is a symmetric-valued tensor measure, it follows that the $\mathscr H_{N-1}$ density of the jump part of $Dw = D^k u$ is of the form 
\[
[w] \otimes \nu_{w}= \lambda_w \otimes (\nu_w)^{\otimes^{k}}
\]
for some $\mathscr H_{N-1}$-integrable map $\lambda_w : J_w \to E$. Using that $|\nu_w|^2 = 1$, we can characterize the $E$-coordinate coefficient $\lambda_w$ directly in terms of $[w]$ and the contraction pairing defined above. Indeed, $\dpr{[w],\nu_w}_{k-1}\cdot |\nu_w|^2 =  \dpr{[w]\otimes \nu_w,\nu_w}_{k}
 = \lambda_w |\nu_\omega|^{2k} = \lambda_w$. 
%As a direct consequence, we deduce that it is possible to recover $[w]$ from the contraction as
%\[
%	[w] = \dpr{[w],\nu_w}_{k-1} \otimes (\nu_w)^{\otimes^{k-1}}.
%\]
In particular, it follows that 
\begin{equation}\label{eq:difference_rep}
	[w] = \lambda_w \otimes  (\nu_w)^{\otimes^{k-1}} \qquad \text{$\mathscr H_{N-1}$ almost everywhere on $J_w$.}
\end{equation}
For consistency, we set $\nabla^0 u \coloneqq u$ so that, when $k=1$, we simply get
\[
	[w] =  \lambda_w = [u]  \qquad \text{$\mathscr H_{N-1}$ almost everywhere on $J_u$.}
\]

With these considerations in mind, we are now in the position to state the following equivalence between Lebesgue-Lusin properties for $k$th order gradients and arbitrary homogeneous operators of order $k$:

	\begin{proposition}\label{prop}
	The following statements are equivalent:
	\begin{enumerate}[label=(\alph*), ref=\alph*,labelsep=1em,left=0pt,itemsep=1em]
	\item \label{a} Let $\Acal$ be a $k$\textsuperscript{th} order linear operator as in~\eqref{eq:A} and let $f : \Omega \to \Fbf_{\Acal}$ be an integrable map. Then
	there exists $u \in BV^k(\Omega;E)$, satisfying
	\begin{gather*}
	\Acal u \; = \; f\,\mathscr L_N \; + \;  \Abb^k(\nu_w) \dpr{[w],\nu_w}_{k-1}   \, \mathscr H_{N-1} \mres J_w\,,
	\end{gather*}
	where $w \coloneqq \nabla^{k-1} u$.  Moreover, $w \in SBV(\Omega;E \otimes \R^{k-1}_\sym)$ and
	\[
	\|u\|_{W^{k-1,1}(\Omega)} + \int_\Omega |\nabla w| \, dx +  \int_{J_w} |[w]| \, d\mathscr H_{N-1} \le C \int_\Omega |f| \, dx 
	\]
 for some constant $C$ depending on $\Acal$.  \label{a}

\item \label{b}	Let $f : \Omega \to  \R^{N^k}_\sym$ be an integrable map. Then
	there exists  $u \in BV^k(\Omega)$ satisfying
	\begin{gather*}
	D^k u \; = \; f\,\mathscr L_N \; + \; [w] \otimes \nu_w  \; \mathscr H_{N-1} \, {\mres \, J_w\,,} 
	\end{gather*}
	where $w \coloneqq \nabla^{k-1} u$. Moreover, $w \in SBV(\Omega;\R^{N^{k-1}}_\sym)$ and
	\[
	\|u\|_{W^{k-1,1}(\Omega)} + \int_\Omega |\nabla w| \, dx + \int_{J_w} |[w]| \, d\mathscr H_{N-1} \le C \int_\Omega |f| \, dx 
	\] 
	for some constant $C$ depending on $N$ and $k$. 
	\end{enumerate}
	\end{proposition}

\begin{proof}[Proof of Proposition~\ref{prop}] 
The implication \eqref{a} $\Longrightarrow$ \eqref{b} is straightforward from the following observations. Firstly, the principal symbol of the $k$th order Hessian is given by $D^k(\xi) = \xi^{\otimes^k}$. 
Therefore $D^k$ has the form~\eqref{eq:A} with $E = \R$ and $F = \R^{N^k}_\sym$, by setting $A_\alpha = M_\alpha$, where the family
$\set{M_\alpha}{|\alpha| =k}$
is the canonical orthonormal basis of $\R^{N^k}_\sym$. %\ado{Moreover, since $u \in W^{k-1,1}_\loc(\Omega)$ and $D^k u$ is symmetric-valued, it follows that the $\mathscr H_{N-1}$ density of the jump part of $D^k u$ is given by (setting $w \coloneqq \nabla^{k-1} u$) $[w] \otimes \nu_{w}= \lambda (\nu_w)^{\otimes^{k}}$ for some Borel map $\lambda : J_w \to \R$. In particular,
%\[
%	[w] \otimes \nu_w = \lambda (\nu_w)^{\otimes^k}.
%\]
%Furthermore, using that $|\nu_w|^2 = 1$, we get 
%\begin{align*}
%	%(w^+ - w^-) \bullet^{k-1} \nu_w & = 
%	\dpr{[w],\nu_w}_{k-1}\cdot |\nu_w|^2 =  \dpr{[w]\otimes \nu_w,\nu_w}_{k}
% = \lambda |\nu_\omega|^{2k} = \lambda.
%\end{align*}
%As a direct consequence, we deduce that 
%\[
%	D^k(\nu_w)\dpr{[w],\nu_w}_{k-1} = \lambda(\nu_w)^{\otimes^k} = [w]  \otimes \nu_w.
%\]
Hence, {from~\eqref{eq:difference_rep} and}~\eqref{a} we conclude that $u \in BV^k(\Omega)$ satisfies {
\[
	D^k u = f\, \mathscr L_N + [w] \otimes \nu_w\, \mathscr H_{N-1} \mres J_w,
\]
with $w = \nabla^{k-1} u \in SBV(\Omega;\R^{N^{k-1}}_\sym)$, satisfying
\[
	\|u\|_{W^{k-1,1}(\Omega)} + \int_\Omega |\nabla w| \, dx + \int_{J_w} |[w]| \le C(D^k) \int_\Omega |f| \, dx .
\]
This proves the first implication since $D^k$ implicitly fixes the spatial dimension $N$.}

We now show  \eqref{b} $\Longrightarrow$ \eqref{a}. The first step is to use the alternative {jet expression} $A[D^k u] = \Acal u$, where $A : E \otimes \R^{N^k}_\sym \to \Fbf_\Acal$ is the (unique) linear map satisfying 
	\begin{equation}\label{eq:1}
		A[e \otimes \xi^{\otimes^k}] = \Abb^k(\xi)[e] \qquad \text{for all $e \in E$\,.}
	\end{equation}
	The existence of $A$ is a direct consequence of the universal property of the tensor product and the $k$-linearity of the principal symbol on the frequency variable $\xi$. 	For the sake of simplicity, let us write $\Ebf_k =E \otimes \R^{N^k}_\sym$. By construction, we have
	\[
		\im A = \spn \set{\Abb^k(\xi)[e]}{\xi \in \R^N, e \in E} = \Fbf_{\Acal} \, .
	\]
Let us consider the Moore-Penrose quasi-inverse $A^\dagger$ associated with $A$. This is an element of $\Hom(\Fbf_\Acal;\Ebf_k)$ satisfying the fundamental identity \begin{equation}\label{eq:2}
	A\circ A^\dagger = \mathbf 1_{\im A} = \mathbf 1_{\Fbf_\Acal} \, .
\end{equation} In particular, by our assumption on $f$, $A^\dagger f(x)$ is well-defined almost everywhere on $\Omega$ and $A^\dagger f \in L^1(\Omega;\Ebf_k)$. 
Now, we make use of  the assumption \eqref{b} over each $E$-coordinate  to find $u \in BV^k(\Omega;E)$ satisfying
	\[
		D^k u  = A^\dagger f  \, \mathscr L_N + D^j w\,, 
	\]
{where $w = \nabla^{k-1} u\, \in SBV(\Omega;\R^{N^k}_\sym)$ satisfies the estimate}
	\begin{equation}\label{eq:aux}
	{\|u\|_{W^{k-1,1}(\Omega)} + \int_\Omega |\nabla w| \, dx + \int_{J_w} |[w]| \, d\mathscr H_{N-1} \le C \int_\Omega |A^\dagger f| \, dx} 
	\end{equation}
for some constant $C$ depending on $N$, $k$ and $\dim \Fbf_\Acal$. 
Here, as usual {$D^j w = \tilde g\, \mathscr H_{N-1} \mres J_w$ is the jump part of $Dw$. Notice that in this case $\tilde g = \mathbf 1_{J_w}  [w] \otimes \nu_w $.} 
	Pre-composing this expression with $A$, we get
	\[
		\Acal u = A [D^k u] = A\left[A^\dagger f\, \mathscr L_N +\tilde g \, \mathscr H_{N-1} \right].
	\]
	Since $A$ is a linear map, we may pull in and distribute $A$ into the densities that belong to the sum of the right-hand side, hence concluding that
	\begin{align*}
	\Acal u =\phantom{:} A\circ A^\dagger f\, \mathscr L_N + A\tilde g \, \mathscr H_{N-1}\stackrel{\eqref{eq:2}}{=}	f\, \mathscr L_N + g\,\mathscr H_{N-1} \, ,
	\end{align*}
where $g =\mathbf 1_{J_w} \cdot  A([w] \otimes \nu_w)$. 
	In passing to the last equality in the formula above, we have used the almost everywhere pointwise restriction $f(x) \in \Fbf_\Acal$. Now, we use the fact that at jump points (cf.~\eqref{eq:difference_rep}), it holds 
	\[
		[w] = \dpr{[w] , \nu_w}_{k-1} \otimes (\nu_w)^{\otimes^{k-1}} \, .
	\] 
From this we obtain that $g =\mathbf 1_{J_w}\cdot \mathbb A^k(\nu_w)\dpr{[w], \nu_w}_{k-1}$ as desired.
Lastly, the bound on the total variation of $D^ku$ follows from~\eqref{eq:aux}, the estimates
	\begin{align*} \int |A^\dagger f|  \, dx & \le \|A^\dagger\|_{\Fbf \to \Ebf_k} \int |f| \, dx \, , \\
		\int |g| \, d\mathscr H_{N-1} & \le \|A\|_{\Ebf_k \to \Fbf} \int |[w]| \, d\mathscr H_{N-1} \, ,
	\end{align*}
and the fact that $\|A^\dagger\|_{\Fbf_\Acal \to \Ebf_k} = \|A^{-1}\|_{\im A \to \Fbf_\Acal}$ depends solely on $\Acal$.\footnote{This equality of norms follows from~\eqref{eq:2}, which implies that the restriction of $A^\dagger$ to $\im A$ is the inverse of the isomorphism $A : (\ker A)^\perp \to \im A$.} This finishes the proof. \end{proof}

\begin{proof}[Proof of Theorems~\ref{thm:1} and~\ref{thm:2}] The proof of Theorem~\ref{thm:1} follows directly from the validity of the Lusin property for gradients (\cite[Thm. 3]{Alberti}),
the previous proposition (with $k =1$),
 and the equivalence in Remark~\ref{rem:coefficients}. The proof of Theorem~\ref{thm:2}, on the other hand, follows from the previous proposition  (with $k = 2$) and~\cite[Theorem 1.4]{Fonseca}. 
\end{proof}

\subsection*{An application for normal currents} Let $m$ be a non-negative integer. An \emph{$m$-dimensional current} $T$ on $\Omega$ is an element of the  continuous dual $\mathscr D_m(\Omega) = \mathscr D^m(\Omega)^*$, where 
$$\mathscr D^m(\Omega) \coloneqq C_c^\infty(\Omega;\wedge^m \R^N)$$
is the space of smooth and compactly supported $m$-forms on $\Omega$. If $m \ge 1$, the (distributional) boundary operator on currents is defined by duality through the rule
\[
	\partial T \in \mathscr D_{m-1}(\Omega)\,, \quad \partial T (\varphi) = T(d\varphi) \; \text{whenever}  \;  \varphi \in \mathscr D^{m-1}(\Omega)\,.
\]
Here, $d$ is the exterior derivative operator acting on $\mathscr D^m(\Omega)$. The boundary operator then maps $\mathscr D_m(\Omega)$ into $\mathscr D_{m-1}(\Omega)$. The mass of a current $T \in \mathscr D_m(\Omega)$ is defined as
\[
	\Mbf(T) = \sup\set{T(\varphi)}{\varphi \in \mathscr D^m(\Omega), \; \sup_{x \in \Omega} |\varphi(x)| \le 1}.
\]
A current $T$ is called \emph{normal} if $T$ is representable by a finite Radon measure and either  $\partial T$ is representable by a finite Radon measure or $m = 0$. In particular, every normal $k$-current is a distribution represented by a measure taking values on the space $\wedge_m \R^N = (\wedge^m \R^N)^*$ of $m$-vectors.
The space of normal $m$-dimensional currents on $\Omega$ is denoted by $\Nbf_m(\Omega)$.

Given an element $\xi \in \R^N \cong \wedge_1 \R^N$, we write $\xi^*$ to denote the $1$-covector
\[
	\xi^*(v) \coloneqq \xi \cdot v, \qquad v \in \R^N.
\] 
The interior multiplication operator $\llcorner: \wedge_{p} \R^N \times \wedge^q \R^N \to \wedge_{p-q} \R^N$ is defined as the adjoint of the exterior multiplication (see e.g.,~\cite[Chapter 1.5]{Federer}):
\[
	\dpr{v \llcorner \alpha,\beta} = \dpr{v,\alpha \wedge \beta} \qquad v \in  \wedge_{p} \R^N, \; \alpha \wedge^q \R^N, \; \beta \in \wedge^{p-q} \R^N.
\]
With this notation in mind, we get the following direct application of Theorem~\ref{thm:1} for the boundary operator acting on normal currents:

\begin{corollary}\label{cor:currents}
Let $m \in [0,N)$ be an integer and let $S : \Omega \to \wedge_m \R^N$ be an integrable $m$-vector field. Then, there exists a normal current $T \in \Nbf_{m+1}(\Omega) \cap SBV(\Omega;\wedge_{m+1}\R^N)$ satisfying
\[
	\partial T \; = \; S \, \mathscr L_N \; + \;  [T] \llcorner \nu_T^* \;  \mathscr H_{N-1} \mres J_T
\]
Moreover,
\[
	\Mbf(T) \, + \, \Mbf(\partial T) \le C_{m,N} \|S\|_{L^1(\Omega)}.
\]
\end{corollary}

\begin{remark}
A similar result holds for integrable $m$-forms with $m \in (0,N]$. More precisely, if $\omega \in L^1(\Omega;\wedge^m \R^N)$, then there exists an $(m-1)$-form $\phi \in SBV(\Omega;\wedge^{m-1} \R^N)$ satisfying (see also~\cite[Proposition 2.1]{Moo})
\[
	d\phi = \omega \, \mathcal L_N + \nu_\phi^* \wedge [\phi] \, \mathscr H_{N-1} \mres J_\phi.
\]
and
\[
	\int_\Omega |u| \, dx + \int_{J_\phi} |[u]| \, d\mathscr H_{N-1} \le C_{m,N} \|\omega\|_{L^1(\Omega)}\,.
\]
\end{remark}

\begin{proof} Since the symbol of the exterior derivative is precisely the exterior multiplication (i.e., $d_m(\xi)\alpha =  \xi^*\wedge \alpha$),  the symbol of the boundary operator $\partial$ on $({m+1})$-vectors is precisely (see, e.g.,~\cite[\S 4.1.7.]{Federer})

\[
	\partial_{m+1} (\xi)e =  -(e \llcorner \xi^*), \qquad \xi \in \R^N, e \in {\wedge_{m+1}} \R^N.
\] 
In particular, the boundary operator $\partial$ defines a  constant-coefficient first-order operator from $\mathscr D_{m+1}(\Omega)$ to $\mathscr D_m(\Omega)$. Since (cf.~\cite[\S 1.5.2]{Federer})
\[
	\im \partial_{m+1}(\xi) =\ker \partial_{m}(\xi) = \spn \{ v_1 \wedge \cdots \wedge v_{m} \; | \; v_1,\dots,v_m \in \xi^\perp\}\,,
\] 
it follows that 
$\wedge_m \R^N = \spn\{ \im \partial_{m+1}(\xi) :|\xi| = 1\} =  \Fbf_\partial$.
We may thus apply Theorem~1 to find $T \in SBV(\Omega;\wedge_m \R^N)$ satisfying the desired properties: The fact that $T \in \Nbf_{m+1}(\Omega)$ follows directly from the fact that $\Mbf(T) = \|T\|_{L^1(\Omega)}$ and that $\Mbf(\partial T) \le |DT|(\Omega)$.
\end{proof}

\begin{corollary}\label{cor:boundaryless} Let $m \in [{1},N)$ be an integer and let $S : \Omega \to \wedge_m \R^N$ be an integrable $m$-vector field. There exists a countably $\mathscr H_{N-1}$  rectifiable set $J$ with (oriented) normal $\nu$, and a Borel $(m+1)$-vector field $g : J \to \wedge_{m+1} \R^N$ satisfying
\[
	\partial (S \, \mathscr L_N - g \llcorner \nu^* \; \mathscr H_{N-1} \mres J) = 0
\]
and 
\[
	\int_{J} |g| \, \mathscr H_{N-1} \le C \int_\Omega |S| \, dx.
\]
\end{corollary}
\begin{proof}
It is sufficient to use the expression for $\partial T$ in the corollary above observe and observe that $\partial(\partial T) = 0$.
\end{proof}

The canonical isomorphism $\iota : \R^N \to \wedge_1 \R^N$, identifying vector fields in $\R^N$ with $1$-vectors, induces an isometry between vector fields with bounded measure divergence and normal $1$-currents. %Namely, it induces the distributional identity 
%\[
%	\diverg f = \sigma \quad \Leftrightarrow \quad \partial (\iota f) = - \sigma.
%\]   
In particular, an exciting and direct application of the previous corollary is the
following \emph{rectifiable completion} for systems of vector fields to systems of solenoidal measures.

\begin{corollary} If $N \ge 2$, then
every vector field ${\vec f} \in L^1(\Omega;\R^{N})$ is the absolutely continuous part of a solenoidal field measure up to a rectifiable measure. More precisely, there exists a countably $\mathscr H_{N-1}$ rectifiable set {$\Gamma \subset \Omega$} and a Borel vector field {$\vec a : \Gamma \to \R^{N}$} satisfying the differential constraint
\[
	\diverg (\vec f\, \mathscr L_N \, + \, \vec a \, \mathscr H_{N-1} \mres \Gamma) = 0
\]
and also the tangential constraint
\[
\vec a(x) \in  \Tan(\Gamma,x) \quad  \text{for $\mathscr H_{N-1}$ almost every $x \in \Gamma$}.
\]
Moreover,  
\[
	  \int_\Gamma |\vec a| \, d\mathscr H_{N-1} \le C \|\vec f\|_{L^1(\Omega)} 
\]
for a constant $C$ that depends only on $N$.
\end{corollary} 
\begin{proof} %If $N = 1$, then the divergence operator is equivalent to the derivative operator $\Acal = d/dt$ and hence the result follows from Theorem~\ref{thm:1} applied to $\vec f \in L^1(\Omega)$. If $N \ge 2$, we argue as follows.
Applying the assertion of Corollary~\ref{cor:boundaryless} (here we are using that $N \ge 2$) to the integrable $1$-vector field {$\iota \vec f$} yields the existence of a countably $\mathscr H_{N-1}$ rectifiable set $\Gamma \subset \Omega$ and an integrable $2$-vector $g : \Gamma \to \wedge_2 \R^N$ satisfying 
\[
	\partial({\iota\vec f} \, \mathscr L_N - g \llcorner \nu^* \, \mathscr H_{N-1} \mres \Gamma) = 0
\]
and
\begin{equation}\label{eq:boundss}
	\int_\Gamma |g| \, \mathscr H_{N-1} \le C \int_\Omega |f|\, dx.
\end{equation}
By the identification discussed before, we conclude that
\[
	\diverg(\vec f \, \mathscr L_N + \vec a \, \mathscr H_{N-1} \mres \Gamma) = 0,
\]
where $\vec a \coloneqq \iota^{-1}(- g \llcorner \nu^*) :\Gamma \to \R^N$ is integrable on $\Gamma$ and satisfies the asserted $L^1$ bounds on $\Gamma$. To see the tangential properties of the $\vec a$, it suffices to observe that $\vec a(x) \cdot \nu(x) = \dpr{g(x) \llcorner \nu(x)^*,\nu(x)^*} = \dpr{g(x), \nu(x)^* \wedge \nu(x)^*} = 0$
for $\mathscr H_{N-1}$ almost every $x \in \Gamma$ (on Lebesgue points of $g$ and $\nu$).
\end{proof}
%\begin{proof}
%
%\end{proof}
%
%
%\begin{proof} 
%Consider the first-order operator 
%\[
%	\curl^* V = \partial_j V_{ij} \qquad V : \Omega \to \R^{N^2}_{\Skew} \, ,
%\]
%whose associated symbol is given by
%\[
%\curl^* (\xi)[P] = P\xi \qquad V\in \R^{N^2}_{\Skew} \, .
%\]
%This defines an operator from $\R^{N^2}_{\Skew}$ to $\R^N$. Notice that $\Fbf_{\curl^*} = \R^N$ so that the essential pointwise restriction is trivial in this case. By Theorem~\ref{thm:1} we may find $V \in BV_\loc(\Omega;\R^{N^2}_{\Skew})$ satisfying
%\[
%	\mu = \curl^* V = f \, \mathscr L_N \; + \; [V^+ - V^-] \nu_V \, \mathscr H_{N-1} \mres J_V
%\]
%and
%\[
%	|DV|(\Omega) \le C \|f\|_{L^1(\Omega)} \, .
%\]
%for some constant $C$ depending only on $N$ (since the symbol of the divergence operator on $\R^N$ depends only on the dimension $N$).
%In particular, 
%Since $V$ is a skew-symmetric matrix field, it follows that
%\[
%	[V^+ - V^-] \nu_V \cdot \nu_V = 0 \quad \Longrightarrow \quad \frac{\mu}{|\mu^s|}  \in \Tan(J_V,x) \; \mathscr H_{N-1} \; \text{almost everywhere\,.}
%\] 
%Moreover, 
%\[
%	\diverg \mu = \diverg \circ \curl^* V = \sum_{i,j = 1}^N\partial_{ij} V_{ij} = 0\,.
%\]
%Setting $R = J_V$ and $g = \mathbf 1_{J_V} \cdot [V^+ - V^-] \nu_V$ yields the sought assertions. 
%\end{proof}

\subsection*{Acknowledgments} I would like to thank Guido De Philippis and Giacomo del Nin for very timely insight and comments about the problem discussed in this note. I would also like to thank the referee for their careful corrections and valuable suggestions. Firstly, for bringing to my attention the work of Fonseca, Leoni, and Paroni on the Lusin property for functions with bounded Hessians (crucial in establishing Theorem~2). And secondly, for suggesting the identification of divergence-bounded fields and normal $1$-currents to give a simpler proof of Corollary~6.


\begin{bibdiv}
    \begin{biblist}


\bib{Alberti}{article}{
   author={Alberti, Giovanni},
   title={A Lusin type theorem for gradients},
   journal={J. Funct. Anal.},
   volume={100},
   date={1991},
   number={1},
   pages={110--118},
   issn={0022-1236},
   %review={\MR{1124295}},
   doi={10.1016/0022-1236(91)90104-D},
}

\bib{rank_one}{article}{
   author={Alberti, Giovanni},
   title={Rank one property for derivatives of functions with bounded
   variation},
   journal={Proc. Roy. Soc. Edinburgh Sect. A},
   volume={123},
   date={1993},
   number={2},
   pages={239--274},
   issn={0308-2105},
  %review={\MR{1215412}},
   doi={10.1017/S030821050002566X},
}

\bib{Alberti2}{article}{
   author={Alberti, Giovanni},
   title={Integral representation of local functionals},
   journal={Ann. Mat. Pura Appl. (4)},
   volume={165},
   date={1993},
   pages={49--86},
   issn={0003-4622},
   %review={\MR{1271411}},
   doi={10.1007/BF01765841},
}

\bib{AD}{article}{
   author={Ambrosio, Luigi},
   author={Coscia, Alessandra},
   author={Dal Maso, Gianni},
   title={Fine properties of functions with bounded deformation},
   journal={Arch. Rational Mech. Anal.},
   volume={139},
   date={1997},
   number={3},
   pages={201--238},
   issn={0003-9527},
    doi={10.1007/s002050050051},
   %review={\MR{1480240}},
}

\bib{APF}{book}{
   author={Ambrosio, Luigi},
   author={Fusco, Nicola},
   author={Pallara, Diego},
   title={Functions of bounded variation and free discontinuity problems},
   series={Oxford Mathematical Monographs},
   publisher={The Clarendon Press, Oxford University Press, New York},
   date={2000},
   pages={xviii+434},
   isbn={0-19-850245-1},
%   review={\MR{1857292}},
}

\bib{Slicing}{article}{
     author={Arroyo-Rabasa, Adolfo},
  title={Slicing and fine properties for functions with bounded {$\mathcal A$}-variation},
  date={2020},
        arxiveprint={
        arxivid={2009.13513},
      },
}

\bib{ADR}{article}{
   author={Arroyo-Rabasa, Adolfo},
   author={De Philippis, Guido},
   author={Rindler, Filip},
   title={Lower semicontinuity and relaxation of linear-growth integral
   functionals under PDE constraints},
   journal={Adv. Calc. Var.},
   volume={13},
   date={2020},
   number={3},
   pages={219--255},
   issn={1864-8258},
   %review={\MR{4116615}},
   doi={10.1515/acv-2017-0003},
}

\bib{AS}{article}{  
  author = {Arroyo-Rabasa, Adolfo},
  author ={Simental, Jos\'e},
  title = {An elementary approach to the homological properties of constant-rank operators},
  journal = {to appear in C. R. Math. Acad. Sci. Paris},
  year = {2023},
  doi={10.5802/crmath.388},
}



\bib{Anna}{article}{
     author={Arroyo-Rabasa, Adolfo},
  author={Skorobogatova, Anna},
  title={A look into some of the fine properties of functions with bounded {$\mathcal A$}-variation},
  date={2019},
        arxiveprint={
        arxivid={1911.08474},
      },
}

\bib{Breit}{article}{
   author={Breit, Dominic},
   author={Diening, Lars},
   author={Gmeineder, Franz},
   title={The Lipschitz truncation of functions of bounded variation},
   journal={Indiana Univ. Math. J.},
   volume={70},
   date={2021},
   number={6},
   pages={2237--2260},
   issn={0022-2518},
%   review={\MR{4359909}},
   doi={10.1512/iumj.2021.70.8742},
}

\bib{Diening}{article}{
   author={Diening, Lars},
   author={Gmeineder, Franz},
   title={Continuity points via Riesz potentials for $\mathbb C$-elliptic
   operators},
   journal={Q. J. Math.},
   volume={71},
   date={2020},
   number={4},
   pages={1201--1218},
   issn={0033-5606},
 doi={10.1093/qmathj/haaa027},
}

\bib{Federer}{book}{
   author={Federer, Herbert},
   title={Geometric measure theory},
   series={Die Grundlehren der mathematischen Wissenschaften, Band 153},
   publisher={Springer-Verlag New York, Inc., New York},
   date={1969},
   pages={xiv+676},
%   review={\MR{0257325}},
}

\bib{Fonseca}{article}{
   author={Fonseca, Irene},
   author={Leoni, Giovanni},
   author={Paroni, Roberto},
   title={On Hessian matrices in the space $BH$},
   journal={Commun. Contemp. Math.},
   volume={7},
   date={2005},
   number={4},
   pages={401--420},
   issn={0219-1997},
%   review={\MR{2166659}},
   doi={10.1142/S0219199705001805},
}

\bib{Greg}{article}{
   author={Francos, Greg},
   title={The Lusin theorem for higher-order derivatives},
   journal={Michigan Math. J.},
   volume={61},
   date={2012},
   number={3},
   pages={507--516},
   issn={0026-2285},
%   review={\MR{2975258}},
   doi={10.1307/mmj/1347040255},
}

\bib{DR}{article}{
   author={De Philippis, Guido},
   author={Rindler, Filip},
   title={On the structure of $\Acal$-free measures and applications},
   journal={Ann. of Math. (2)},
   volume={184},
   date={2016},
   number={3},
   pages={1017--1039},
   issn={0003-486X},
%   review={\MR{3549629}},
   doi={10.4007/annals.2016.184.3.10},
}



\bib{Murat}{article}{
   author={Murat, Fran\c{c}ois},
   title={Compacit\'{e} par compensation: condition n\'{e}cessaire et suffisante de
   continuit\'{e} faible sous une hypoth\`ese de rang constant},
   language={French},
   journal={Ann. Scuola Norm. Sup. Pisa Cl. Sci. (4)},
   volume={8},
   date={1981},
   number={1},
   pages={69--102},
   issn={0391-173X},
%   review={\MR{616901}},
}

\bib{Moo}{article}{
   author={Moonens, Laurent},
   author={Pfeffer, Washek F.},
   title={The multidimensional Lusin theorem},
   journal={J. Math. Anal. Appl.},
   volume={339},
   date={2008},
   number={1},
   pages={746--752},
   issn={0022-247X},
%   review={\MR{2370691}},
   doi={10.1016/j.jmaa.2007.07.041},
}

\bib{Raita}{article}{
   author={Rai\c{t}\u{a}, Bogdan},
   title={Potentials for $\Acal$-quasiconvexity},
   journal={Calc. Var. Partial Differential Equations},
   volume={58},
   date={2019},
   number={3},
   pages={Paper No. 105, 16},
   issn={0944-2669},
%   review={\MR{3958799}},
   doi={10.1007/s00526-019-1544-x},
}





    \end{biblist}
    \end{bibdiv}
\end{document}